\documentclass[12pt, a4paper,twoside]{article}
\usepackage{a4,amssymb,amsfonts,verbatim,xypic,times}
\usepackage{yfonts}
\usepackage[pagewise]{lineno}
\usepackage[english]{babel}
\usepackage{amssymb,latexsym}
\usepackage{amsmath}
\usepackage{latexsym}

\newcommand{\M}{\mathbb{M}}
\newcommand{\uS}{{\not\hspace{-0.1cm}\Sigma}}
\newcommand{\ug}{{\not\hspace{-0.1cm}\gamma}}
\newcommand{\un}{{\not\hspace{-0.1cm}\nabla}}

\newcommand{\uD}{{\not\hspace{-0.1cm}D}}

\newcommand{\Imm}{{\mathfrak{Imm}}}
\newcommand{\CC}{{\mathbb C}}

\newcommand{\RR}{{\mathbb R}}

\newcommand{\Spec}{{\mathrm Spec\ }}

\newcommand{\Ekt}{\mathbb{E}(\kappa,\tau)}

\newcommand{\cZ}{{\mathcal Z}}

\newcommand{\Ric}{{\mathrm {Ric}}}

\renewcommand{\phi}{\varphi}

\newcommand{\Spin}{\mathrm{Spin}}

\newcommand{\End}{\mathrm{End}}

\newcommand{\tr}{\mathrm{tr}}
\newcommand{\<}{\left\langle}       
\renewcommand{\>}{\right\rangle}

\newcommand{\Spinc}{\mathrm{Spin^c}}

\newcommand{\id}{\mathrm{Id}} 
\newcommand{\vol}{\mathrm{vol}}

\newtheorem{example}{Examples}[section]
\newtheorem{thm}{Theorem}[section]
\newtheorem{lemma}[thm]{Lemma}
\newtheorem{prop}[thm]{Proposition}
\newtheorem{cor}[thm]{Corollary}
\newtheorem{remark}[thm]{Remark}
\newtheorem{remarks}[thm]{Remarks}
\newtheorem{definition}[thm]{Definition}
\newtheorem{notation}[thm]{Notation}
\newtheorem{exabout:ample}[thm]{Example}

\begin{document}
\title{{\bfseries Hypersurfaces immersed in special Spin$^c$ manifolds by first eigenspinors}}
\author{Roger Nakad}


\maketitle
\begin{center}
Department of Mathematics and Statistics,\\
Faculty of Natural and Applied Sciences,\\
Notre Dame University-Louaize\\
P.O. Box 72, Zouk Mikael, \\
Zouk Mosbeh, \\
Lebanon
\end{center}
\begin{center}
 {\bf rnakad@ndu.edu.lb}
\end{center}
\vskip 0.5cm
\begin{center}
 {\bf Abstract}
\end{center}

Let $M$ be a closed orientable hypersurface of dimension $n$, with nonwhere vanishing mean curvature $H$, immersed into a Riemannian $\Spinc$ manifold $\cZ$ carrying a parallel spinor field.  The first eigenvalue $\lambda_1(\uD)$ (with the least absolute value) of the induced Dirac operator $\uD$ of $M$ satisfies the $\Spinc$ B\"{a}r inequality \cite{Ba3, nakadthese, NR1},
\begin{eqnarray*}
\lambda_1^2 (\uD) \le  \frac{n^2}{4 \ \vol(M)}\int_M H^2 dV,
\end{eqnarray*}
where $\vol(M)$ is the volume of $M$ and $dV$ is the volume form 
of the manifold $M$. In this paper, we classify hypersurfaces $M$ that satisfy the equality case in the $\Spinc$ B\"{a}r inequality when $\cZ = (0,+\infty) \times P$ is the cone over a Riemannian $\Spinc$ manifold $P$ carrying a real Killing spinor, under two conditions: one being a Ricci condition on $\cZ$, and the second one the curvature of the auxiliary line bundle associated with the $\Spinc$ structure on $\cZ$. More precisely, we prove that $M$ are the slices $\{s\} \times P$, where $s \in (0,+\infty)$. In the special case, when $\cZ=\mathbb R^{n+1}$, i.e., the cone over the sphere, which is a Spin manifold with a parallel spinor, the classification result was previously obtained by Hijazi and Montiel \cite{HM2012}.



{\bf Acknowledgement}: This work was completed at the International Centre for Theoretical Physics (ICTP) in Trieste, Italy, where the author spent two weeks as part of the “ICTP-INdAM Collaborative Grants and Research in Pairs – 2025” program. The author gratefully acknowledges the financial support of ICTP and sincerely thanks the institute for its continuous support and hospitality during the research stay.

\section{Introduction}

In recent years, researchers in Submanifold Theory and experts in Riemannian and Spin Geometry have increasingly studied Dirac operators naturally associated with the extrinsic geometry of hypersurfaces, exploring their spectral properties to gain deeper insights into the geometry and topology of these hypersurfaces (see \cite{AF, Am1, Am2, Ba3, BFLPP, FLPP, HMR2, HMZ1, HMZ2, HMZ3, HM2012, HM2014, T} and related references). These Spin Geometry techniques have led to significant advances in extrinsic geometry, including the study of constant mean curvature (CMC) or minimal surfaces in homogeneous $3$-spaces from Thurston's classification of 3-dimensional geometries, the Willmore conjecture, the Willmore functional and Alexandrov-type theorems. Remarkably, Spin Geometry tools \textemdash especially special spinor fields and the Dirac operator\textemdash have played a crucial role in many of these developments, inspiring further research directions.

When shifting from classical Spin geometry to $\Spinc$ geometry, the framework becomes more general, giving rise to a range of challenges and obstacles. This is because the $\Spinc$ structure depends not only on the geometry of the manifold but also on the connection (and hence the curvature) of the auxiliary line bundle associated with the chosen $\Spinc$ structure. From a physical point of view, spinors model fermions while $\Spinc$ spinors can be
interpreted as fermions coupled to an electromagnetic field. Extending the use of Spin methods in submanifold studies to the $\Spinc$ setting allows for the exploration of a broader range of ambient geometric structures. Indeed, every Spin manifold is trivially $\Spinc$ and manifolds such as CR, K\"{a}hler, complex, and Sasaki manifolds are naturally $\Spinc$, often with special spinor fields. For instance, from an extrinsic viewpoint, O. Hijazi, S. Montiel, and F. Urbano constructed $\Spinc$ structures on K\"{a}hler-Einstein manifolds with positive scalar curvature, carrying K\"{a}hlerian Killing spinors \cite{HMU}. When restricted to minimal Lagrangian submanifolds, these spinors impose topological and geometric constraints on the submanifolds. Further applications of $\Spinc$ geometry in extrinsic geometry can be found in \cite{nakadthese, NR}. We will now highlight an example of such results that is particularly relevant to the problem we will address in this paper.

{\bf Problem Setting:}  Let $M$ be an $n$-dimensional closed oriented hypersurface, isometrically immersed in a $(n+1)$-dimensional Riemannian $\Spinc$ manifold $\cZ$, as is typical for many ambient spaces in Submanifold Theory. Under these conditions, $M$ inherits an induced
 $\Spinc$ structure. We will denote the spinor bundle on $M$ by $\uS M$ and the Dirac operator by $\uD$. This situation occurs, for example, when the ambient space $\cZ$ is $\mathbb{CP}^{n+1}$, which is always a $\Spinc$ manifold\textemdash Spin when $n+1$ is odd and non-Spin when $n+1$ is even. In fact, it carries a canonical $\Spinc$ structure endowed with a parallel spinor. More generally, a complete, simply connected manifold admits a parallel $\Spinc$ spinor if and only if it is isometric to the Riemannian product of a simply connected K\"{a}hler manifold (with its canonical or anti-canonical $\Spinc$ structure) and a simply connected Spin manifold with a parallel spinor \cite{19, HM}. Examples of such Riemannian Spin manifolds admitting non-trivial parallel spinor fields include Euclidean spaces, Calabi-Yau manifolds, hyper-K\"{a}hler manifolds, certain $7$- and $8$-dimensional special Riemannian manifolds, and their Riemannian products \cite{Wa1}.

Let $\psi$ be a parallel spinor field on $\cZ$, which we can normalize to have unit constant length. We also denote by $\psi$ the restriction of this parallel spinor to the hypersurface $M$. Similar to the approach used by Bleecker and Weiner for the Laplacian, B\"{a}r \cite{Ba3} used this spinor field $\psi$ as a test spinor in the variational characterization (Rayleigh quotient) of the lowest eigenvalues of the Spin Dirac operator on $M$, deriving an upper bound for the first eigenvalue of the Spin Dirac operator $\uD$. This inequality, known as the B\"{a}r inequality, was subsequently extended to the eigenvalues of the hypersurface $\Spinc$ Dirac operator $\uD$ on $M$ \cite{nakadthese, NR1}. More precisely, we have the following result:
\begin{thm}[\cite{Ba3, nakadthese, NR1}.
]\label{barSpinc}
Let $M$ be a compact orientable hypersurface immersed into a Riemannian $\Spinc$ manifold $\cZ$ carrying a parallel spinor. 
Let $\lambda_1(\uD)$ be the first eigenvalue (with the least absolute value) of the
induced Dirac operator $\uD$ of $M$. Then

\begin{eqnarray}
\lambda_1^2 (\uD) \le  \frac{n^2}{4 \ \vol(M)}\int_M H^2 dV,
\label{inequalitybar}
\end{eqnarray}
where $\vol(M)$ is the volume of $M$, $dV$ is the volume form 
of the manifold $M$ and $H$ the mean curvature of $M$.
\end{thm}

\begin{center}
 {\bf Problem Statement:} {\it Can we characterize hypersurfaces $M$ that satisfy the equality case in the B\"{a}r Inequality (\ref{inequalitybar}) ?}
 \end{center}
For several years, it was an open problem to determine when equality holds in the B\"{a}r inequality (\ref{inequalitybar}), even for the Spin Dirac operator. Indeed, equality directly implies that the test spinor field $\psi$ must be an eigenspinor of $\uD$. It is evident that this condition implies that the mean curvature $H$ is constant. However, this information alone is insufficient unless the ambient manifold is $\mathbb R^{n+1}$ and the hypersurface $M$ is embedded. In this case, we can apply the well-known Alexandrov theorem \cite{Al}, which states: {\it Any closed hypersurface embedded in Euclidean space with constant mean curvature must be a hypersphere}. As a result of this theorem, during the 1960s and 1970s, many geometers studying submanifolds believed that spheres might be the only hypersurfaces immersed in Euclidean space with constant mean curvature (see \cite{Ho}). However, since 1983, Hsiang et al. \cite{HTY} have constructed many examples of such hypersurfaces that are not round spheres, particularly in even dimensions of the ambient space, while Wente \cite{We} and Kapouleas \cite{Ka1, Ka2} have provided examples in three-dimensional space.

In summary, if equality holds in (\ref{inequalitybar}), the mean curvature $H$ of the hypersurface $M$ is constant. Consequently, as previously noted, the restriction of the parallel spinor field $\psi$ becomes an eigenspinor of the Dirac operator $\uD$, corresponding to one of the eigenvalues $\pm \frac n2 H$ and this eigenvalue must be the one with the least absolute value. Recently, hypersurfaces $M$ that achieve equality in the B\"{a}r inequality have been classified when $\cZ=\mathbb R^{n+1}$ \cite{HM2012}, revealing that $M$ must be a hypersphere. A careful examination of the proofs indicates that the use of the position vector field $x$ defined on $\mathbb R^{n+1}$ is crucial. In other words, a key ingredient is the existence of a vector field $x \in T(\mathbb R^{n+1})$ such that $\nabla_Xx = X$ for all $X \in  T(\mathbb R^{n+1})$. But the existence of such a vector field is a characteristic feature of the cone over an arbitrary $n$-dimensional Riemannian manifold $P$, that is the $(n + 1)$-dimensional Riemannian manifold $\cZ$ defined as $\cZ =(0, +\infty) \times P$ endowed with the following warped product metric:
$$dr\otimes dr + r^2 \<., .\>_P,\ \ \ \ \forall r \in (0, +\infty).$$
In fact, the vector field $x$ given by
$$x = r \frac{\partial}{\partial r}$$
plays the role of the vector position field in the Euclidean space, that is $\nabla_Xx = X$ for all
$X \in T \cZ$. Note that for $P = \mathbb S^n$, we have $\cZ= \mathbb R^{n+1}$, and this is the only case where
the singularity at the vertex of the cone can be removed.

B\"{a}r made an insightful observation (see \cite{Ba1}) that the cone $\cZ$ is a Riemannian Spin manifold with a non-trivial parallel spinor if the factor $P$ is a Riemannian Spin manifold equipped with a non-trivial real Killing spinor. Assuming for simplicity that $P$ is simply connected, and aside from the basic case $P= \mathbb S^n$, $P$ could be either an Einstein-Sasaki Riemannian manifold, a $7$-dimensional Spin(7)-manifold, or a $6$-dimensional nearly-K\"{a}hler manifold of constant type 1 (see \cite{Ba1}). For any of these choices of $P$, the corresponding cone $\cZ$ is a Riemannian Spin manifold endowed with a non-trivial parallel spinor and a position vector field. In the $\Spinc$ case, a similar result holds: the cone $\cZ$ is a Riemannian $\Spinc$ manifold with a non-trivial parallel spinor if $P$ is taken to be a Riemannian $\Spinc$  manifold endowed with a non-trivial real Killing spinor \cite{Mo}. However, despite the existence of a position vector field, and as pointed previously, several challenges arise since the $\Spinc$ structure depends not only on the geometry of the manifold but also on the connection (and thus the curvature) of the auxiliary line bundle associated with the $\Spinc$ structure. For instance, Spin manifolds with real Killing spinors (resp. parallel spinors) are Einstein (resp. Ricci flat) but $\Spinc$ manifolds with real Killing spinors (resp. parallel spinors) are not necessarily Einstein (resp. Ricci flat). We will answer the previously stated problem in the $\Spinc$ context, and indeed, we will prove the following theorem:


\begin{thm}\label{thmmain2}
Let $P^n$ be a closed Riemannian $\Spinc$ manifold of dimension  $n$ carrying a Killing spinor. Consider the cone $\cZ = (0, +\infty)\times P$ over $P$ and assume that  
\begin{eqnarray}\label{co1}
\mathrm{R}^\cZ \geq c_{n+1} \vert\Omega^\cZ\vert,
\end{eqnarray}
\begin{eqnarray}\label{co2}
\mathrm{Ric}^\cZ (X, X) \geq 0, \ \ \ \ \text{for all}\ \ \ \ X \in \Gamma(T\cZ)
\end{eqnarray}
where $\mathrm{Ric}^\cZ$ is the Ricci of $\cZ$, $\mathrm{R}^\cZ$ denotes the scalar curvature of $\cZ$,  $i\Omega^\cZ$ the curvature $2$-form of the auxiliary line bundle associated with the $\Spinc$ structure ($\Omega^\cZ$ is a real $2$-form on $M$), and $c_{n+1} = 2[\frac {n+1}{2}]^{\frac 12}$. Here $[\frac {n+1}{2}]$ is the integer part of $\frac {n+1}{2}$. Let $M$ be a closed hypersurface immersed into $\cZ$  satisfying the equality case of the $\Spinc$ B\"ar type inequality (\ref{inequalitybar}).  Then, $M$ is a slice $\{s\} \times P$, where $s \in (0, \infty)$. We shall say that the only compact constant mean curvature hypersurfaces of $\cZ = (0, +\infty)\times P$ immersed by first 
eigenspinors are slices $\{s\} \times P$.
\end{thm}

If $P$ is a $\Spin$ manifold with a real Killing spinor then condition (\ref{co1}) is always satisfied because the manifold is Einstein with positive scalar curvature. Moreover, the cone carries a parallel
spinor so it is Ricci flat, and the scalar curvature is zero. Hence, the condition $0 = \mathrm{R}^\cZ \geq c_{n+1}\vert \Omega\vert^\cZ= 0$ is also satisfied. 
Again, in the very special case is when $P = \mathbb S^n$, and thus $\cZ = \mathbb R^{n+1}$, the result was proved by Hijazi and Montiel in \cite{HM2012}. Two non-$\Spin$ examples are provided in Section \ref{nonspin}. In fact, we prove that the only compact  hypersurfaces of $\cZ = (0, +\infty) \times \mathbb B^3$, for some Berger Spheres $\mathbb B^3$, that satisfy the equality case in (\ref{inequalitybar}) are slices $\{ s \} \times \mathbb B^3$. Berger Spheres are $\Spinc$ manifolds with a real Killing spinor. They are also Spin but without real Killing spinors. Also, we construct explicitly a $5$-dimensional non-$\Spin$, non-Einstein Sasakian manifold such that its cone $\cZ$ fullfills the two conditions of Theorem \ref{thmmain2}. 

As a consequence of the Theorem \ref{thmmain2}, we can prove the following Alexandrov-type theorem:
\begin{thm} [{\bf Alexandrov-type Theorem}] Let $P$ be a compact Riemannian $\Spinc$ manifold of dimension  $n$ carrying a real Killing spinor.
Consider $\cZ = (0, +\infty) \times P$ the cone over $P$ and assume that  conditions  (\ref{co1}) and (\ref{co2}) are satisfied. Then, the only compact embedded hypersurfaces  into $\cZ$  of constant positive  mean curvature are  slices $\{s\} \times P$.
\label{thmmain3}\end{thm}


\section{Preliminaries}
In this section, we briefly recall basic notions about the Dirac operator on $\Spinc$ manifolds and their hypersurfaces. More details can be found in \cite{BHMM}, \cite{Fr}, \cite{LM}, \cite{mon}, \cite{nakadthese} and \cite{IEspinc}.

{\bf The Dirac operator on Spin$^c$ manifolds.} Let $(\cZ^{n+1}, g= \<., .\>)$ be a Riemannian $\Spinc$ manifold of dimension $n+1$ without boundary. On such a manifold, we have  a Hermitian complex vector bundle $\Sigma \cZ$ endowed with a natural scalar product (denoted also by $\<., .\>$) and a connection $\nabla $ that parallelizes the metric. That is,
\begin{eqnarray}\label{eq1}
  \label{eq1} X\<\psi, \phi\> = \<\nabla_X\psi, \phi\> + \<\psi, \nabla_X\phi\>,
  \end{eqnarray} 
   for any $\psi, \phi \in \Gamma(\Sigma \cZ)$ and $X\in \Gamma(T\cZ)$. This complex vector bundle, called the $\Spinc$ bundle, is endowed with a Clifford multiplication denoted by ''$\gamma$", $\gamma\colon T\cZ \rightarrow \mathrm{End}_{\mathbb C} (\Sigma \cZ)$, such that at every point $x \in \cZ$, it defines an irreducible representation of the corresponding Clifford algebra. Hence, the complex rank of $\Sigma \cZ$ is $2^{[\frac {n+1}{2}]}$. The Hermitian scalar product $\<., .\>$ is compatible with this Clifford multiplication  and the spinorial connection $\nabla$. That is,
 \begin{eqnarray}
  \label{eq2} \<\gamma(X)\psi, \gamma(X)\psi\> = \vert X\vert^2 \<\psi, \phi\>,
  \end{eqnarray}  
  \begin{eqnarray}
  \label{eq3} \nabla_X(\gamma(Y)\psi) = \gamma(\nabla_XY) \psi + \gamma(Y)\nabla_X\psi,
  \end{eqnarray} 
for any $\psi, \phi \in \Gamma(\Sigma \cZ)$ and $X, Y \in \Gamma(T\cZ)$. 

Given a $\Spinc$ structure on $(\cZ^{n+1}, g)$, one can prove that the determinant line bundle $\mathrm{det} (\Sigma \cZ)$ has a root of index $2^{[\frac{n+1}{2}]-1}$. We denote
by $L$ this root line bundle over $\cZ$ and it is called the auxiliary line bundle associated with the $\Spinc$ structure. Locally, a $\Spin$ structure always exists. We denote by $\Sigma' \cZ$ the (possibly globally non-existent)
spinor bundle. Moreover, the square root of the auxiliary line bundle $L$
always exists locally. But, $\Sigma\cZ = \Sigma' \cZ \otimes {L}^{\frac 12}$ exists globally.  This essentially means that, while
the spinor bundle and ${L}^{\frac 12}$
may not exist globally, their tensor product (the $\Spinc$  bundle) is
defined globally. Thus, the connection $\nabla$ on $\Sigma \cZ $ is the twisted connection of the one on the
spinor bundle (coming from the Levi-Civita connection) and a fixed connection on $L$. 

We may now define the Dirac operator $D$ acting on the space of smooth
sections of $\Sigma\cZ$  by the composition of the metric connection and
the Clifford multiplication. In local coordinates this reads as
\begin{eqnarray*} 
D = \sum_{j=1}^{n+1} \gamma(e_j) \nabla_{e_j},
 \end{eqnarray*}
where $\{e_1,\ldots,e_{n+1}\}$ is a local oriented orthonormal tangent frame. It is a first order elliptic operator, formally self-adjoint with respect to the $L^2$-scalar product and satisfies the Schr\"odinger-Lichnerowicz formula 
\begin{eqnarray}
D^2 = \nabla^*\nabla + \frac 14 \mathrm{R}^{\cZ} \id_{\Gamma (\Sigma \cZ)}+ \frac{i}{2}\gamma(\Omega^\cZ),
\label{sl}
\end{eqnarray}
where $\nabla^*$ is the adjoint of $\nabla$ with respect to the
 $L^2$-scalar product and $\gamma( \Omega^{\cZ}) $ is the extension of the Clifford multiplication to differential forms. 
By definition, the
function $\<D\psi, \psi\>$ is nothing but the trace of a symmetric two-tensor T$^\psi$,
associated with the spinor field $\psi$, the so-called energy-momentum
tensor  which is defined as follows
\begin{eqnarray}\label{EMT}
T^\psi (X, Y) = \frac 12 \<\gamma(X)\nabla_Y\psi + \gamma(Y)\nabla_X\psi, \psi\>,
\end{eqnarray} 
for any $X, Y \in \Gamma(\cZ)$. 

Every almost complex manifold $(\cZ^{2m = n+1}, g, J)$ 
of complex dimension $m$ has a canonical $\Spinc$ structure. In fact, the complexified cotangent bundle 
$T^*\cZ\otimes \CC = \Lambda^{1,0} \cZ \oplus \Lambda^{0,1}\cZ$ 
decomposes into the $\pm i$-eigenbundles of the complex linear extension of the complex structure. Thus, the spinor bundle of the canonical $\Spinc$ structure is given by $$\Sigma \cZ = \Lambda^{0,*} \cZ =\oplus_{r=0}^m \Lambda^{0,r}\cZ,$$
where $\Lambda^{0,r}\cZ = \Lambda^r(\Lambda^{0,1}\cZ)$ is the bundle of $r$-forms of type $(0, 1)$. The auxiliary
 line bundle of this canonical $\Spinc$ structure is given by  $L = (K_\cZ)^{-1}= \Lambda^m (\Lambda^{0,1}\cZ)$, 
where $K_\cZ$ is the canonical bundle of $\cZ$ \cite{Fr, mon, nakadthese}. Let $\ltimes$ 
be the K\"{a}hler form defined by the complex structure $J$, i.e. $\ltimes (X, Y)= g(X, JY)$ for all vector 
fields $X,Y\in \Gamma(T\cZ).$ The auxiliary line bundle $L= (K_\cZ)^{-1}$ has a canonical holomorphic 
connection induced from the Levi-Civita connection whose curvature form is given 
by $i\Omega^\cZ = i\rho^\cZ$, where $\rho$ is the Ricci $2$-form given by $\rho(X, Y) = \mathrm{Ric} (X, JY)$. For the canonical $\Spinc$ structure and when $\cZ$ is a  K\"{a}hler manifold, admits parallel spinors (constant functions) 
lying in $\Lambda^{0,0}\cZ$\cite{19}. Of course, we can define another $\Spinc$ structure for which the spinor bundle is 
given by $\Lambda^{*, 0} \cZ =\oplus_{r=0}^m \Lambda^r (T_{1, 0}^* \cZ)$ and the auxiliary line bundle by $K_\cZ$. 
This $\Spinc$ structure will be called the anti-canonical $\Spinc$ structure.

{\bf Spin$^c$ hypersurfaces and the Gauss formula.} Here, we follow the same notations as in \cite{HM2012}. 
We compare the restriction $\uS M$ of the $\Spinc$ bundle $\Sigma\cZ$
of the $\Spinc$ manifold $\cZ$ to an orientable hypersurface $M$ immersed
into $M$ and its Dirac-type operator $\uD$ to the intrinsic spinor bundle
$\Sigma M$ of the induced $\Spinc$ structure of $M$ and its Dirac operator $D_M$. 
Let $\un$ be the Levi-Civita connection associated with the induced
Riemannian metric on $M$. The Gauss formula says that
\begin{eqnarray}\label{GF}
\un_X Y = \nabla_XY  - \<AX, Y\> N,
\end{eqnarray}
where $X, Y$ are vector fields tangent to the hypersurface $M$, the vector
field $N$ is a global unit field normal to $M$ and $A$ is for the shape
operator corresponding to $N$, that is, for all $X\in \Gamma(TM)$,
\begin{eqnarray}\label{sff}
\nabla_X N = -AX.
\end{eqnarray}
We have that the restriction $\uS M : = \Sigma \cZ_{\vert_M}$ 
is a left module over $\CC l(M)$ for the induced Clifford multiplication
$\ug : \CC l(M) \longrightarrow \End(\uS M)$ given by 
\begin{eqnarray}\label{cmh}
\ug (X) \psi = \gamma(X)\gamma(N)\psi,
\end{eqnarray}
for every $\psi\in \Gamma(\uS M)$ and $X\in \Gamma(TM)$. Consider on $\uS M$ the Hermitian
metric $\<\cdot, \cdot\>$ induced from that of $\Sigma\cZ$. This metric immediately satisfies the compatibility condition (\ref{eq2})
 if one puts on $M$ the Riemannian
metric induced from $\cZ$ and the Clifford multiplication $\ug$ defined in
(\ref{cmh}). The connection $\un$ from Equation (\ref{GF}) induces a connection $\un$ on $\uS M$ given by
\begin{eqnarray}\label{gaussspinc}
\un_X\psi = \nabla_X\psi - \frac 12 \gamma(AX)\gamma(N)\psi = \nabla_X\psi - \frac 12 \ug(AX)\psi,
\end{eqnarray}
for every   $\psi\in \Gamma(\uS M)$ and $X\in \gamma(TM)$. Note that a suitable use of the
Gauss formula (\ref{GF}) shows that the compability conditions (\ref{eq1}) and (\ref{eq3})
are also satisfied for $(\uS M, \ug, <\cdot, \cdot>, \un)$. Denote by $\uD : \Gamma(\uS M) \longrightarrow \Gamma(\uS M)$ the Dirac operator associated with
the Dirac bundle $\uS M$ over the hypersurface, which is defined by
$$\uD\psi = \sum_{j=1}^n \ug (e_j) \un_{e_j}\psi,$$
where   $\psi \in \Gamma(\uS M)$ and $\{e_1,\cdots, e_n\}$ is a local orthonormal frame of
$TM$. Since, as we can see from (\ref{cmh}) and (\ref{gaussspinc}), the connection $\un$ does
not depend on the choice of the normal $N$ on $M$ while $\ug$ changes sign
for the opposite orientation, hence the same is true for $\uD$ .
It is a well known fact that $\uD$ is a first order elliptic differential
operator which is formally $L^2$-self-adjoint. By (\ref{gaussspinc}), for any spinor field
  $\psi\in \Gamma(\uS M)$, we have
\begin{eqnarray}\label{gaussdirac}
\uD \psi = \frac n2 H\psi - \gamma(N) \sum_{j=1}^n \gamma(e_j)\nabla_{e_j}\psi,
\end{eqnarray}
where $H = \frac 1n \tr (A)$ is the mean curvature of $M$ corresponding to the
orientation $N$. Note that if the spinor field $\psi$  is the restriction to the
hypersurface $M$ of a spinor field on the ambient space $\cZ$, then both
spinor fields in $\Gamma(\Sigma \cZ)$ and $\Gamma(\uS M)$, will be denoted by the same symbol. 
For any spinor $\psi \in \Gamma(\uS M)$   and any tangent vector field $X\in \Gamma(TM)$, 
the following relations hold
\begin{lemma}[\cite{HM2012}. Lemma 5]\label{lemma5} For any $X \in \Gamma(TM)$ and $\psi \in \Gamma(\uS M)$, we have
$$\un_X(\gamma(N)\psi) = \gamma(N)\un_X\psi,$$
$$\uD(\gamma(N)\psi) = -\gamma(N) \uD\psi.$$
\end{lemma}
In particular, the second of these equalities implies that, when M is a compact manifold
without boundary, the spectrum of $\uD$ must be symmetric with respect
to zero. This does not occur necessarily for $D_M$. When the dimension $n$ of $M$ is even we have $(\uS M, \ug, \uD) \equiv (\Sigma M, \gamma_M, D_M)$
 and the decomposition $\uS M = \uS M^+ \oplus \uS M^-$, given
by $\uS M^\pm := \{ \eta \in \uS M, i\gamma(N)\eta = \pm \eta\}$ corresponds, up to the
above identification, to the chirality decomposition of the spinor
bundle $\Sigma M$. Hence $\uD$ interchanges $\uS M^+$ and $\uS M^-$, as one can
also easily see from Lemma \ref{lemma5}. When $n$ is odd, the chirality decomposition of $\Sigma\cZ$ into positive
and negative spinors induces an orthogonal and 
$\ug$, $\uD$-invariant
decomposition $\uS M = \uS M_+ \oplus \uS M_-$, with $\uS M_{\pm} := (\Sigma^\pm \cZ)_{\vert_M}$, in
such a way that $(\uS M_\pm, \ug, \uD_{\vert_{{{\not\hspace{-0.00011cm}S}M}_{\pm}}}) \equiv (\Sigma M, \pm\gamma_M, \pm D_M)$. 
 Moreover, we have the following two isomorphisms: $\gamma(N): \uS M_{\pm} \longrightarrow \uS M_\mp$. Furthermore, if $M$ is compact without boundary, then the discrete spectra
of the elliptic self-adjoint operators $\uD$ and $D_M$, satisfy the following:
\begin{enumerate}
\item $\Spec \uD$ is symmetric with respect to zero.
\item If $n$ is even, we have $\Spec \uD = \Spec D_M$ 
\item If $n$ is odd, we have  
$\Spec \uD = \Spec D_M \cup (-\Spec D_M)$. 
\end{enumerate}
The above considerations regarding the relation between $\uD$ and $D_M$ allow us to work henceforth
 only with $\uS M$and $\uD$, 
and shall refer
to them as the spinor bundle and the Dirac operator of the induced $\Spinc$ structure of the hypersurface $M$.

\begin{remark}
 Although, by definition, sections of the restricted bundle
$\uS M$ are defined only on the hypersurface 
$M$ and not on the whole
space $\cZ$, we may consider derivatives of the form $\nabla_X \psi$, where $\nabla$ is the
spinorial Levi-Civita connection on the ambient space and where  $\psi\in \Gamma(\uS M)$, provided that the field $X$
 is tangent to $M$. In fact, at any point
of $M$, the spinor field   can be locally extended to $\cZ$, we compute the
covariant derivative of the extension and restrict again to $M$. The result
does not depend on the chosen local extension. Alternatively, if we
think of the restricted bundle $\uS M$ as the spinor bundle $\Sigma M$, according
to the above identifications, via the Gauss formula (\ref{gaussspinc}), we may define
the connection $\nabla$ in terms of the spinorial Levi-Civita connection $\un$,
the Clifford multiplication  $\ug$ on $\Sigma M$ and the shape operator $A$ of the
hypersurface.
\end{remark}
Applying (\ref{gaussdirac}), we obtain the counterpart of Proposition 7 in \cite{HM2012}, this time within the $\Spinc$ framework:
\begin{prop}\label{byfirsteigenspinors}
Let $M$ be an orientable hypersurface immersed into a Riemannian $\Spinc$ manifold $\cZ$ admitting a parallel
spinor field $\psi \in \Gamma(\Sigma\cZ) $. Then the restriction   $\psi\in \Gamma(\uS M)$ satisfies
$$ \uD \psi = \frac n2 H \psi.$$
If $H$ is constant, then $\frac n2 H$ belongs to the spectrum of $\uD$ , 
hence $\lambda_1(\uD ) \leq \frac n2 \vert H \vert$. If $\lambda_1(\uD ) = \frac n2 \vert H \vert$, 
we shall say that {\bf $M$ is immersed in $\cZ$ by first eigenspinors}.\end{prop}
\begin{remark}\label{RR}
Consider the case of the round unit hypersphere $M = \mathbb S^n$
as an embedded hypersurface of $\mathbb R^{n+1}$, an ambient space which admits
a maximal number of independent parallel spinor fields. It is well known
that the eigenvalues of the Dirac operator $\uD$ of the unique $\Spin$
structure on $\mathbb S^n$ with the lowest absolute value are $\pm \frac n2$ (see \cite{Gi}), 
both with multiplicity $2^{[\frac{n+1}{2}]}$. Moreover it is also clear that
$H = \pm 1$ according to the choice of the orientation, hence in this case $\lambda_1(\uD) = \frac n2 \vert H \vert$.
\end{remark}

\section {Totally umbillic hypersurfaces of $\Spinc$ manifolds carrying parallel spinors}
In this section, we relate the Ricci curvature of $\cZ$ to the curvature  of the auxiliary line bundle $L$, when $M$ be a totally umbillic hypersurface of $\cZ$ of constant mean curvature $H$. 
\begin{prop}\label{hyper_totally_umb}
 Let $\cZ^{n+1}$ be a $\Spinc$-manifold carrying a parallel spinor $\psi$. Let $M$ be a totally umbillic hypersurface of $\cZ$ of dimension $n$ and of constant mean curvature $H$. Then $\Ric^\cZ (N, e_i) = 0$ on $M$ and  $\Vert \mathrm{Ric}^\cZ (N, N) \Vert =  \Vert N \lrcorner \Omega^\cZ \Vert$ on $M$.
\end{prop}
{\bf Proof.} Let $\phi= \psi|_M$.
 By $A=H g$ and \eqref{gaussspinc} we get $0=\nabla_X\phi +\frac{H}{2}\gamma(X)\gamma(N)\phi$  for all $X\in \Gamma(TM)$. Moreover, we have $\uD^2\phi=\frac{n^2}{4}H^2\phi$.
 Together with $ \Delta^M\phi = n\frac{H^2}{4}\phi$
 and the Lichnerowicz formula for hypersurfaces (see \cite{GrNa}, Lemma 3.1), we obtain
\begin{eqnarray}\label{Lich_M_hyper} 
n(n-1)\left(\frac{H^2}{4}\right)\phi-\frac{\mathrm{R}^M}{4}\phi= \frac{i}{2} \ug(\Omega^M)\phi,
\end{eqnarray}
where $\mathrm{R}^M$ denotes the scalar curvature of $M$. The Schr\"odinger formula (\ref{sl}) gives that $$0=\frac{\mathrm{R}^\cZ}{4}\phi+\frac{i}{2}\gamma({\Omega^\cZ})\phi.$$
Using the relation between $\Omega^\cZ$ and $\Omega^M$ (see \cite{IEspinc}) given by
 $$\frac{i}{2}\gamma({\Omega^\cZ})\phi = \frac{i}{2}\ug({\Omega^M})\phi - \frac{i}{2} \ug(N \lrcorner {\Omega^\cZ})\phi,$$  one has 
\begin{eqnarray}\label{26}
0=\frac{\mathrm R^\cZ}{4}\phi+\frac{i}{2} \ug(\Omega^M)\phi-\frac{i}{2} \ug(N\lrcorner\Omega^\cZ)\phi.
\end{eqnarray}
Substracting \eqref{Lich_M_hyper} from Equation (\ref{26}), we get 
\[ -n(n-1)\frac{H^2}{4}\phi + \frac{\mathrm R^M}{4} \phi=-\frac{i}{2} \ug(N\lrcorner\Omega^\cZ)\phi + \frac{\mathrm R^\cZ}{4} \phi\]
Using that  $\mathrm R^M =  \mathrm R^\cZ - 2 \mathrm{Ric}^\cZ (N, N) + n^2 H^2 - \vert II \vert^2,$ yields
\begin{eqnarray}\label{Id2}
0=\frac{\Ric^\cZ(N,N)}{2}\phi-\frac{i}{2} \ug(N\lrcorner\Omega^\cZ)\phi.
\end{eqnarray}
We recall the Ricci identity on $\cZ$ \cite{GrNa}:
\begin{eqnarray}\label{RC}
\sum_{k=1}^{n+1} \gamma(e_k){\mathrm R}(e_k, X)\Psi = \frac 12 \gamma(\mathrm{Ric}^\cZ
(X) )\Psi - \frac i2 \gamma(X\lrcorner\Omega^{\cZ})\Psi,
\end{eqnarray}
for any spinor field $\Psi$, $X \in\Gamma(T\cZ)$, and where $\mathrm{Ric^\cZ}$ denotes the Ricci curvature of $\cZ$ and $\mathrm{R}$ is the spinorial curvature. Applying (\ref{RC}) to the parallel spinor $\psi$ and the vector $N$, we get
\begin{eqnarray*}
\mathrm{Ric}^\cZ (N, N) \gamma(N)  \psi + \sum_{j=1}^n \mathrm{Ric}^\cZ (N, e_j) \gamma(e_j)\psi -i\gamma(N \lrcorner  \Omega^\cZ)\psi = 0
\end{eqnarray*}
Taking the Clifford multiplication of the last identity with $N$ gives 
\begin{eqnarray}\label{Id3}
\mathrm{Ric}^\cZ (N, N)\psi + \sum_{j=1}^n \mathrm{Ric}^\cZ (N, e_j) \ug(e_j)\psi -i \ug(N \lrcorner \Omega^\cZ)\psi = 0
\end{eqnarray}
Now (\ref{Id3}) and (\ref{Id2}) give that $\Ric^\cZ (N, e_i) = 0$ on $M$ and $\Vert \mathrm{Ric}^\cZ (N, N) \Vert^2 =  \Vert N \lrcorner \Omega^\cZ \Vert^2$ on $M$.
\begin{remark}
The condition of $M$ being of constant mean curvature can be removed in Proposition \ref{hyper_totally_umb} if we assume that the dimension of $\cZ$ satisfies $n+1\geq 5$ (see  \cite{GrNa} for more details). In fact, every connected totally umbilical hypersurface  of a Riemann-
ian $\Spinc$ manifold of dimension $n+1   \geq 5$ carrying a real (including parallel) Killing
spinor is of constant mean curvature. However, there exist totally umbilical connected hypersurfaces with non-constant mean curvature in Riemannian $\Spinc$ manifolds of dimension $3$ or $4$ and
carrying parallel spinors \cite{GrNa}.
\end{remark}

In \cite{Mo}, Montiel provided a characterization of totally umbilical hypersurfaces with constant mean curvature in certain special manifolds. This characterization presents, a priori, three possibilities, two of which must be excluded when the mean curvature is nonzero. In this case, the result \textemdash which we will use to characterize hypersurfaces immersed in $\cZ$ via first eigenspinors\textemdash can be stated as follows:

\begin{thm} \cite[Lemma 4]{Mo}\label{Montiel-classification}
Let $M^n$ be a closed orientable hypersurface immersed in a totally umbilical way and with nonzero constant mean curvature into a Riemannian manifold $\cZ^{n+1}$ endowed with a non trivial closed conformal vector field $V$ such that the Ricci curvature satisfies 
\begin{eqnarray}\label{CO}
\vert V \vert^2 \Ric^\cZ(X, X) - \vert X \vert^2 \Ric^\cZ(V, V) \geq 0,
\end{eqnarray}  
for all $X \in \Gamma(T \cZ)$. If a unit vector field normal to $M$ attains the equality case, then $M^n$ is a leaf of the foliation.
\end{thm}


\section{Comformal change of the hypersurface Dirac operator and the energy-momentum tensor}\label{conformalgeometry}

Suppose now that the mean curvature $H$ of the immersed (connected)
hypersurface $M$ of the $\Spinc$ manifold $\cZ$ has no zeros. Fix
the orientation and the corresponding unit normal field so that $H > 0$. 
Consider on $M$ the conformal metric  $\<\cdot, \cdot\>^H = H^2 \<\cdot, \cdot\>$. If we denote
by $\Sigma^H M$ the $\Spinc$ bundle corresponding to this conformal metric and
the same $\Spinc$ structure. We shall identify the two spinor bundles $\Sigma M$
and $\Sigma^H M$  using the isomorphism given in \cite{Hi1}. Due to the existence
of this isometric identification, from now on, we shall denote the two
spinor bundles by the same symbol $\Sigma M$. With this identification, the
corresponding Clifford multiplications and $\Spinc$ connections are related
as follows:
\begin{eqnarray}\label{conf-covariant}
{\ug}^H = H\ug, \ \ \un_X^H - \un_X = -\frac{1}{2H} \ug
(X)\ug(\nabla H) + \<X, \nabla H\>,
\end{eqnarray}
for all $X\in \Gamma(TM)$. From this, we can easily find the following relation
between the two Dirac operators $\uD$ and $\uD^H$ on $\Sigma M$ relative to the two
conformal metrics (see \cite{Hi1})
\begin{eqnarray}\label{conf-dirac}
\uD^H (H^{-\frac{n-1}{2}} \phi) = H^{-\frac{n+1}{2} }\uD\phi,
\end{eqnarray}
for any $\phi \in \Gamma(\uS M)$. 
\begin{remark}\label{Hconst}
When $H$ is constant, the metrics $<\cdot, \cdot>^H$ and $<\cdot, \cdot>$ are homothetic. In this case, the conformal
 covariance (\ref{conf-dirac}) becomes the following change of scale: $\uD^H = \frac{1}{H}\uD$.
\end{remark}
Suppose now that there exists a non-trivial parallel spinor field  $\psi \in \Gamma(\Sigma \cZ)$ and take 
 $\phi$ in (\ref{conf-dirac}) its restriction   $\psi \in \Gamma(\uS M)$. Using Proposition
 \ref{byfirsteigenspinors}, we obtain
$$\uD^H (H^{-\frac{n-1}{2}}\psi) = H^{-\frac{n+1}{2}} \uD\psi = \frac n2 H^{-\frac{n-1}{2}} \psi.$$
From this equality we can deduce the following result

\begin{prop}\label{barconf}
Let $M$ be an orientable hypersurface immersed into
a Riemannian $\Spinc$ manifold $\cZ$ admitting a parallel spinor field $\psi \in \Gamma(\Sigma \cZ)$ 
and suppose that the mean curvature of M, after a suitable
choice of the unit normal, satisfies $H > 0$. Denote by  $\psi^H$ the spinor field $H^{-\frac{n-1}{2}} \psi \in \Gamma(\uS M)$. 
Then we have $\uD^H \psi ^H = \frac{n}{2}\psi^H$. As a consequence, $\frac{n}{2}$ belongs to the spectrum of $\uD^H$, hence 
\begin{eqnarray}\label{upperH}
\lambda_1(\uD^H) \leq \frac n2.
\end{eqnarray}
\end{prop}
\begin{remark}\label{remark4}
It is important to note that, when the mean curvature H
of the hypersurface is constant and $\lambda_1 (\uD)=  \frac n2 H$, that is, when $M$
 is immersed by first eigenspinors (see Proposition \ref{byfirsteigenspinors}), then from Remark \ref{Hconst}, we get
$\lambda_1(\uD^H) = \frac n2$. That is, in this case, the equality is achieved in
(\ref{upperH}).
\end{remark}
\begin{remark}
 Proposition \ref{barconf} can be stated in terms of a spectral functional
$\mathfrak F_1$ defined on the space of immersions $\mathfrak i: M \longrightarrow \cZ$ of a compact
n-dimensional manifold $M$ into a Riemannian $\Spinc$ $(n+1)$-dimensional
manifold $(\cZ, <\cdot, \cdot>)$ with non vanishing mean curvature, that we shall
denote by $\Imm^+ (M,\cZ)$. This functional is defined by
$$\mathfrak F_1 : \mathfrak i \in \Imm^+ (M, \cZ) \longrightarrow \lambda_1(\uD^{H_{\mathfrak i}}) \in \RR,$$
where $\lambda_1(\uD^{H_{\mathfrak i}})$ is the first non-negative eigenvalue of the Dirac operator
$\uD^{H_{\mathfrak i}}$ of the induced $\Spinc$ structure corresponding to the conformal
metric on $M$ given by
$$\<\cdot, \cdot\>^{H_{\mathfrak i}} = H^2_{\mathfrak i} \mathfrak{i}^* \<\cdot, \cdot\>.$$
This functional $\mathfrak F_1$ is clearly invariant by homotheties of the ambient
space metric. Proposition \ref{barconf} can be paraphrased by saying that this
functional is bounded from above by $\frac n2$ provided that there exists a
non-trivial parallel spinor field on the ambient space $\cZ$. In fact, when
$\cZ = \RR^{n+1}$, from Remark \ref{RR} and Proposition \ref{barconf}, we have that this bound
is achieved by hyperspheres of arbitrary radius, so $\frac n2$  is a maximum. In general, from Remark \ref{remark4}, 
we can see that all constant mean curvature
hypersurfaces in Riemannian $\Spinc$ manifolds with non-trivial parallel
spinors immersed by first eigenspinors attain this maximum, too.

\end{remark}
\section{Hypersurfaces with critical eigenvalue $\frac n2$ of $\uD^H$}

This section closely follows Section 5 in \cite{HM2012}, but we have rewritten it for clarity, emphasizing the key difference: the curvature term of the auxiliary line bundle that arises in the $\Spinc$ setting.
\\\\
Suppose that $\mathfrak i: M^n \longrightarrow \cZ^{n+1}$ is an immersion of a closed orientable
$n$-dimensional manifold $M$, into a Riemannian $\Spinc$ manifold $\cZ$, whose
mean curvature function $H$ is positive and such that $\lambda_k(\uD^H) = \frac n2$, 
for
some $k \in \mathbb N$. Note that here we are not assuming the existence of a non
trivial parallel spinor field on $\cZ$. Consider an analytic variation $\mathfrak i(t)$ of
the immersion $\mathfrak i$, that is, a family of immersions analytically indexed
by $t\in ]-\epsilon, +\epsilon[$ such that $\mathfrak i(0) = \mathfrak i$. 
Since they are homotopic, all the
immersions $\mathfrak i(t)$ of the variation induce on $M$  the same $\Spinc$ structure
from the $\Spinc$ structure on $\cZ$. 
Moreover, the variation induces an
analytic deformation of the Riemannian metric $\<\cdot, \cdot\>^H$ on $M$, namely,
\begin{eqnarray}\label{Ht}
\<\cdot, \cdot\>^H (t) = \<\cdot, \cdot\>(t)^{H(t)} = H(t)^2 \<\cdot, \cdot\> (t),
\end{eqnarray}
where $H(t)$ (which remains to be positive if $\epsilon$ is small enough) and
 $\<\cdot, \cdot\>(t)$ are the mean curvature and the induced metric
of the immersion $\mathfrak i(t)$, respectively. This deformation of the metric determines an
analytic family of Dirac operators $\uD^H(t)$ given by  $\uD^H(t) = \uD(t)^{H(t)}$, 
each acting on the different spinor bundles $\uS^{H(t)}_t M$ corresponding to
the common induced $\Spinc$ structure and the metrics  $\<\cdot, \cdot\>^{H}(t)$, where
we are using the notation introduced in Section \ref{conformalgeometry} for conformal changes
of the metric and where the spinor bundle $\uS_tM$ corresponds to the induced
$\Spinc$ structure and to the metric $ \<\cdot, \cdot\>(t)$. In order to compare
the different deformed Dirac operators 
$\uD^H(t)$  with the Dirac operator $\uD^H(0) = \uD^H$
 associated to the original immersion, it is necessary to
relate the corresponding spinor bundles $\uS_t^{H(t)} M$  with $\uS_0^{H(0)} M = \uS^H M$.
There is no canonical way to do this, but Bourguignon and Gauduchon
determined in \cite{BG} a specific way to define an isometry between two
spinor bundles corresponding to the same spin structure but to two different metrics on a given manifold. This isometry allowed them to
consider the two associated Dirac operators as two operators acting
on the same bundle and to obtain relation between them (Th\'eor\`eme 21 in \cite{BG}). This isometry has been extented to $\Spinc$ manifolds \cite{IEspinc} (Theorem 5.1). 
It is also important to observe that, when the two metrics are conformally related, the Bourguignon and Gauduchon 
identification coincides with that of Section 4. Then, using this isometry,
we may see the one-parameter family $\uD^H (t)$ as an analytic family
of Dirac operators acting on the same spinor bundle $\uS^H M \simeq \uS M$.\\\\
Now, we are interested in the first variation of the eigenvalues $\lambda_k (\uD^H)$
with respect to the analytic deformation $\mathfrak i(t)$ of the immersion $\mathfrak i$, that
is, with respect to the variation of the metric $\<\cdot, \cdot\>^H (t)$ on $M$. In other
words, we want to compute the derivatives
$$
\frac{d}{dt}\vert_{t=0} \lambda_k(t),$$
provided that they exist, where $
\lambda_k(t) = \lambda_k(\uD^H (t)).$
For this, one can see in \cite{BG} that the Rellich-Kato theory of unbounded
self-adjoint operators can be applied to the analytic family
$\uD^H(t)$ acting on $\uS M$. Denote by $m$
 the dimension of the eigenspace $E_k(\uD^H)$ 
associated with the eigenvalue $\lambda_k(\uD^H)$. Then, for each $t \in ]-\epsilon, +\epsilon[$, 
there exist $m$ eigenvalues $\Lambda_{k,1}(t), \cdots, \Lambda_{k,m}(t)$ of $\uD^H (t)$ associated with an 
$L^2 (M, \<\cdot, \cdot\>^H(t))$-orthonormal family of eigenspinors
$\psi_{k,1}(t), \cdots, \psi_{k,m}(t) \in \Gamma(\uS M)$ of $\uD^H(t)$, both $\Lambda_{k,i}(t)$
 and $\psi_{k,i}(t)$ depending 
analytically on $t$, for $i = 1\cdots m$, such that 
$$\Lambda_{k,1} (0) = \cdots =\Lambda_{k,m}(0) = \lambda_k(\uD^H),$$ and
$$\uD^H(t) \psi_{k,i}(t) = \Lambda_{k,i}(t)\psi_{k,i}(t),$$
for all $t \in ]-\epsilon, +\epsilon[$ and $i=1, \ldots, m$.\\
Then, in order to study the critical points, 
we need to know the
derivatives $\Lambda^{'}_{k,i}(0)$ of the $m$ branches $\Lambda_{k,i}(t)$ passing through $\lambda_k(\uD^H)$ at
$t = 0$. The values of these derivatives were computed in \cite{BG, IEspinc}: for each $i = 1,\ldots, m$,
$$
\frac{d}{dt}\vert_{t=0} \Lambda_{k,i}(t) = -\frac 12 \int_M \<T_{\psi_i}^H, h\>^H dV^H.$$
where $dV^H$ 
is the Riemannian density of the metric $\<\cdot, \cdot\>^H$ on $M$ and
 $\psi_i = \psi_{k,i} (0)$ is a fortiori an $L^2(M, \<\cdot, \cdot\>^H)$-orthonormal basis of the
eigenspace $E_k(\uD^H)$. Here $T^\psi_H$
  stands for the energy-momentum tensor
(as defined in (\ref{EMT})) of the spinor $\psi \in \Gamma(\uS M)$, calculated with respect
to the conformal metric $\<\cdot, \cdot\>^H$ on $M$, that is, for all vector fields $X, Y$
tangent to $M$:
\begin{eqnarray}\label{EMTc}
 T^H_\psi (X, Y) = \frac 12  \<\ug^H(X)\un_Y^H \psi + \ug^H(Y)\un_X^H \psi, \psi\>,
\end{eqnarray}
and $h$ is the variational field associated with the metric variation $\<\cdot, \cdot\>^H (t)$ that is
\begin{eqnarray}\label{hfield}
h = \frac{d}{dt}\vert_{t=0} \<\cdot, \cdot\>^H (t).
\end{eqnarray}
It can be shown, see \cite{BG, IEspinc}, that these $m$ derivatives $\frac{d}{dt}\vert_{t=0} \Lambda_{k,i} (t)$
are just the eigenvalues of the quadratic form $Q_h$ defined on the $m$-
dimensional eigenspace $E_k(\uD^H)$ as follows: 
\begin{eqnarray}\label{Qh}
 \phi \in E_k (\uD^H) \longrightarrow Q_h (\phi) = - \frac 12 \int_M \<T^H_\phi, h\>^H dV^H.
\end{eqnarray}

Our next goal is to compute $Q_h(\phi)$ in terms of the metric $\<\cdot, \cdot\>$
induced on $M$. First note that since the functional $\mathfrak F_k: \mathfrak{i}\in \Imm^+ (M,\cZ) \to \lambda_k(\uD^{H_{\mathfrak i}})$ is geometric,
it is invariant under reparametrizations of $M$, so there is no loss of
generality by assuming that the variation $\mathfrak i(t)$ of the immersion $\mathfrak i$ is a
normal variation, that is,
\begin{eqnarray}\label{fN}
 \frac {d}{dt}\vert_{t=0} \mathfrak i(t) = f N,
\end{eqnarray}
where $f$ is a smooth function defined on $M$. \\
Then, using (\ref{Ht}) and (\ref{hfield}), the variational tensor field $h$ corresponding
to the metric deformation $\<\cdot, \cdot\>^H (t)$ is precisely 
$$ h  =2H (\frac{d}{dt}\vert_{t=0} H(t)) \<\cdot, \cdot\> + H^2 \frac{d}{dt}\vert_{t=0} \<\cdot, \cdot\> (t).$$
Now, it is a well-known fact in Submanifold Theory that the variation
of the induced metric associated with the normal variation (\ref{fN}) of the
immersion is given by
$$\frac{d}{dt}\vert_{t=0} \<\cdot, \cdot\> (t) = -2 f \<A \cdot , \cdot\>,$$
and that the corresponding variation of
the mean curvature is given by (see \cite{Mo}):
$$\frac{d}{dt}\vert_{t=0} H(t) = \frac 1n (\bigtriangleup f + \vert A\vert^2 f +\mathrm{Ric}^\cZ (N, N)f).$$
Finally we get
\begin{eqnarray} \label{scal}
h = \frac 2n H( \bigtriangleup f +\vert A\vert^2 f + \mathrm{Ric}^\cZ (N, N)f) \<\cdot, \cdot\> - 2H^2 f \<A\cdot, \cdot\>. 
\end{eqnarray}
We also need to compute the energy-momentum tensor $T^H_\phi$  associated
with an eigenspinor $\phi \in E_k(\uD^H)$. By definition of this tensor (see (\ref{EMT}) and (\ref{EMTc})) 
and the relation between the Levi-Civita connections and
the Clifford multiplications of two conformally related metrics established
in (\ref{conf-covariant}), we get the following simple covariance rule \cite{Hi2}
\begin{eqnarray} \label{TH}
T^H_\phi = HT_\phi.
\end{eqnarray}
Choose any eigenspinor $\phi \in  E_k(\uD^H)$. Since $H > 0$, we have that
$$\phi = \psi^H = H^{-\frac{n-1}{2}}\psi$$  
for a unique  $\psi\in \Gamma(\uS M)$. Then by the conformal covariance (\ref{conf-dirac}) and
the fact that $\phi$ is an eigenspinor of $\uD^H$, we get
$$ \lambda_k(\uD^H) \phi = \uD^H \phi = H^{-\frac{n+1}{2}} \uD\psi,$$
hence
\begin{eqnarray} \label{DH}
\uD \psi = \lambda_k(\uD^H) H\psi = \frac n2 H\psi.
\end{eqnarray}
Similarly, the conformal covariance (\ref{TH}) of the energy-momentum tensor
can be expressed in terms of the spinor field  . It suffices to use
the properties of   and the conformal covariance
$$T_{g\phi} = g^2 T_\phi, \ \ \forall g\in C^\infty (M)$$
which can be deduced from the definition of the tensor $T_\phi$ (see (\ref{EMT})).
So, we have
\begin{eqnarray} \label{THfinal}
T^H_\phi = HT_{H^{-\frac{n-1}{2}}\psi} = H^{-n+2} T_\psi.
\end{eqnarray}

Finally we are in position to compute the integrand in the quadratic
form (\ref{Qh}) which controls the first derivative of the functional $\mathfrak F_k$ at
$t = 0$. First, we observe that
\begin{eqnarray} \label{THscal}
\<T^H_\phi, h\>^H = \frac{1}{H^4} \<T^H_\phi, h\>,
\end{eqnarray}
since this a scalar product of two-tensors. Now, using (\ref{scal}), (\ref{DH}), (\ref{THfinal})
and (\ref{THscal}), we deduce the expression
\begin{eqnarray} \label{finalexpression}
\<T^H_\phi, h\>^H = \frac{1}{H^n} \Big(\bigtriangleup f + \vert A\vert^2 f +\mathrm{Ric}^\cZ(N, N)f) \vert\psi\vert^2 - 2f\<A, T_\psi\> \Big).
\end{eqnarray}
In order to compute the contraction $\<A, T_\psi\>$, take the squared length
in the first equality of (\ref{gaussspinc}), where the two connections $\nabla$ and $\un$ on
the restricted bundle $\uS M$ are related, considering that, for all $X, Y \in
\Gamma(TM)$, the corresponding version of (\ref{EMT}) is now,
$$T_\psi(X, Y) = \frac 12  <\ug(X)\un_Y\psi + \ug(Y)\un_X\psi, \psi>.$$
Hence we get
$$\<A, T_\psi\> = -\vert\nabla\psi\vert^2 + \vert\un\psi\vert^2 + \frac 14 \vert A\vert^2 \vert\psi\vert^2.$$
which when inserted in (\ref{finalexpression}) and recalling that the Riemannian measures
of the two conformal metrics $\<\cdot, \cdot\>$ and $\<\cdot, \cdot\>^H$
$$dV^H = H^n dV,$$
we finally obtain
$$\<T^H_\phi, h\>^H dV^H = \Big( (\bigtriangleup f +\frac 12 \vert A\vert^2 f +\mathrm{Ric}^\cZ(N, N)f) \vert\psi\vert^2 + 2 (\vert\nabla\psi\vert^2 - \vert\un\psi\vert^2)f\Big)dV.$$
Integration by parts of this formula, for $\psi$  belonging to the $m$-dimensional
space $H^{-\frac{n-1}{2}} E_K(\uD^H)$, shows that the quadratic form $Q_h$ (henceforth denoted
by $Q_f$ since it depends only on $f$) takes the form
$$Q_f (\psi ) = -\frac 12 \int_M \Big( \bigtriangleup \vert\psi\vert^2 + \frac 12 \vert A\vert^2 \vert\psi\vert^2 + \mathrm{Ric}^\cZ(N, N) \vert\psi\vert^2 + 2 \vert\nabla\psi\vert^2 - 2 \vert\un\psi\vert^2\Big) fdV,$$
and controls the derivatives $\frac{d}{dt}|_{t=0} \Lambda_{k,i}$,  $i = 1, \ldots, m$ of the $m$ analytic
branches of the eigenvalue $\lambda_k(\uD)$ of $\uD^H(t)$ at $t = 0$. Note that the
Laplacian of the squared length of an eigenspinor is easily computable.

In fact, from the compatibility condition (\ref{eq1}), we immediately have
$$ \bigtriangleup \vert\psi\vert^2 = 2\vert\un\psi\vert^2 + 2 \<\tr \un^2\psi, \psi\>$$
But the rough Laplacian tracer $\un^2$ acting on sections of the $\Spinc$ bundle
$\uS M$ satisfies the celebrated Schr\"{o}dinger-Lichnerowicz formula
$$\tr \un^2 = \frac 14 \mathrm{R}^M + \frac i2 \ug(\Omega^M) - \uD^2,$$
where $\mathrm {R}^M$ is the scalar curvature of $M$ and $\Omega^M$ the curvature of the line bundle on $M$. This and (\ref{DH}) imply
$$ \bigtriangleup\vert\psi\vert^2  = 2\vert\un\psi\vert^2 + \frac 12( \mathrm{R}^M -n^2H^2)\vert\psi\vert^2 + \<i\ug(\Omega^M)\psi, \psi\>.   $$
Finally, the Gauss equation relating the curvature tensor of the ambient
space $\cZ$ and that of the hypersurface $M$ reads
$$\mathrm{R}^M = \mathrm{R}^\cZ -2\mathrm{Ric}^\cZ(N, N) +n^2H^2 - \vert A\vert^2.$$
Hence, the following equality for functions on $M$:
$$\bigtriangleup \vert\psi\vert^2 = 2 \vert\un\psi\vert^2 +
\big( \frac 12 \mathrm{R}^\cZ- \mathrm{Ric}^\cZ(N, N) - \frac 12 \vert A\vert^2\big) \vert\psi\vert^2 + 
\<i\ug(\Omega^M)\psi, \psi\>.$$
Thus, we get 
$$Q_f (\psi) =
 -\int_M (\vert\nabla\psi\vert^2 + \frac 14 R^{\cZ} \vert\psi\vert^2 + \frac i2 \<\ug(\Omega^M)\psi, \psi\>) fdV.$$
Moreover, we have \cite{IEspinc}
$$\gamma(\Omega^\cZ)\psi = \ug (\Omega^M) \psi - \ug(N \lrcorner \Omega^\cZ)\psi.$$
Finally,
\begin{eqnarray*}
Q_f (\psi) =
 -\int_M (\vert\nabla\psi\vert^2 + \frac 14 R^{\cZ} \vert\psi\vert^2 &+& \frac i2 \<\gamma(\Omega^\cZ)\psi, \psi\> \\ &&+
\frac i2 \<\ug(N\lrcorner\Omega^\cZ)\psi, \psi\> ) fdV.
\end{eqnarray*}

With the above facts in mind, we are now ready to give the following
result which allows the characterization of immersed hypersurfaces for
which the eigenvalue $\frac n2$ of $\uD^H$ is a critical point of the functional $\mathfrak F_k$.
\begin{prop}
Let $M$ be a compact orientable hypersurface immersed
into a Riemannian $\Spinc$ manifold $\cZ$ and assume that the mean curvature
$H$ of $M$ is positive. Moreover, suppose that $\frac n2 = \lambda_k (\uD^H)$ is an eigenvalue with multiplicity $m$ of the Dirac operator, corresponding to the conformal
metric $\<\cdot, \cdot\>^H$ and to the induced $\Spinc$ structure.
Consider a normal variation $\mathfrak i(t)$ of the immersion $M$ with variational
field given by $fN$, where $f$ is a smooth function on $M$, and let $\mathfrak F_k(t)$
be the corresponding variation of the $k$-th eigenvalue of $\uD^H$. Then,
the values of the derivatives of $\mathfrak F_k$ lie among the m eigenvalues of the
quadratic form
\begin{eqnarray*}
\psi \in H^{\frac{n-1}{2}} E_k(\uD^H) \mapsto Q_f (\psi) =
 -\int_M (\vert\nabla\psi\vert^2 &+& \frac 14 R^{\cZ} \vert\psi\vert^2 + \frac i2 \<\gamma(\Omega^\cZ)\psi, \psi\> \\ && +
\frac i2 \<\ug(N\lrcorner\Omega^\cZ)\psi, \psi\> ) fdV.
\end{eqnarray*}
\end{prop}
We choose $f\equiv 1$ that is, considering parallel variations
of the hypersurface $M$. If we assume that the operator $\mathrm{R}^\cZ + 2i\gamma(\Omega^\cZ) + 2i \gamma(N) \gamma(N\lrcorner\Omega^\cZ)$ is non-negative, it is then obvious that the quadratic expression depending on the spinor  $\psi$, inside
the parentheses, is non-negative. Thus, the quadratic form $Q_1$ is non-positive on the $m$-dimensional vector space on which it is defined. Then all its eigenvalues
will be non-positive. This leads to the following necessary condition
for the immersion to be critical with respect to the functional $\mathfrak F_k = \lambda_k(\uD^H)$.
\begin{thm}\label{thm10}
 Let M be an $n$-dimensional compact orientable hypersurface
of a Riemannian $\Spinc$ manifold $(\cZ, \<\cdot, \cdot\>)$. Suppose that
$$\mathrm{R}^\cZ + 2i\gamma(\Omega^\cZ) + 2i \gamma(N) \gamma(N\lrcorner\Omega^\cZ)$$
 is non-negative and that the mean curvature
$H$ of $M$ is positive with respect to a suitable choice of the normal. If $\frac n2 = \lambda_k (\uD^H)$ 
belongs to the spectrum of the Dirac operator $\uD^H$ of the
metric $H^2 \<\cdot, \cdot\>$ and it is critical for all the variations of the hypersurface
$M$, then 
$$ \mathrm{R}^\cZ \vert\psi\vert^2+ 2i \<\gamma(\Omega^\cZ)\psi, \psi\> + 2i \<\gamma(N)\gamma(N\lrcorner\Omega^\cZ)\psi, \psi\> =0\ \ \ \ \text{on} \ \ \ \  M,$$ $$\nabla \psi =0 ,$$
for any spinor field  $\psi$ on $M$
satisfying
$$\uD \psi = \frac n2 H\psi.$$
\end{thm}

\begin{remark}
\label{remarkk}
Let's give examples where the operator 
\begin{eqnarray}\label{operator}
\mathrm{R}^\cZ + 2i\gamma(\Omega^\cZ) + 2i \gamma(N) \gamma(N\lrcorner\Omega^\cZ)
\end{eqnarray}
 is non-negative.
\begin{enumerate}
\item  When $\cZ$ is $\Spin$, the operator (\ref{operator}) is equal to the scalar curvature $\mathrm{R}^\cZ$ \cite{HM2012} and hence being non-negative is equivalent to say that $\mathrm{R}^\cZ$ is non-negative.
\item When $N \lrcorner \Omega^{\cZ} = 0$ and $\mathrm{R}^\cZ \geq c_{n+1} \vert \Omega^\cZ \vert $, the operator (\ref{operator}) is non-negative. In fact, using the Cauchy-Schwartz inequality, we  have \cite{HM}
\begin{eqnarray}\label{CSomega}
\<i\gamma(\Omega^\cZ)\Psi, \Psi\> \geq - \frac{c_{n+1}}{2} \vert\Omega^\cZ\vert \vert\Psi\vert^2,
\end{eqnarray}
For any $\Psi \in \Gamma(\Sigma\cZ)$. Thus, it gives
\begin{eqnarray*}
\mathrm{R}^\cZ \vert \psi\vert^2 &+& 2i\<\gamma(\Omega^\cZ)\psi, \psi\> + 2i \<\gamma(N) \gamma(N\lrcorner\Omega^\cZ)\psi, \psi\>\\ &=& \mathrm{R}^\cZ \vert \psi\vert^2+ 2i\<\gamma(\Omega^\cZ) \psi, \psi\>\\ &\geq&  (\mathrm{R}^\cZ - c_{n+1} \vert \Omega^\cZ\vert) \vert \psi\vert^2, 
\end{eqnarray*}
which is non negative.
\end{enumerate}
\end{remark}

\begin{cor}\label{cor11}
Let $M$ be an $n$-dimensional compact orientable hypersurface
of a Riemannian $\Spinc$ manifold $(\cZ, \<\cdot, \cdot\>)$ whose mean curvature
$H$ is positive. Assume that  
 $$\mathrm{R}^\cZ \geq c_{n+1}\vert\Omega^\cZ\vert,\ \ \ N \lrcorner \Omega^{\cZ} = 0,$$ and $\cZ$ carries a parallel spinor.  
Then the first non-negative eigenvalue $\lambda_1(\uD^H)$ of
the metric $H^2 \<\cdot, \cdot\>$  is at most $\frac n2$. If the equality holds, then 
$$\nabla\psi =0\ \ \ \ \text{on}\ \ M,$$
$$\mathrm{R}^\cZ = c_{n+1}\vert\Omega^\cZ\vert \ \ \ \ \text{on}\ \ \ M,$$
for
any spinor field $\psi$   on $M$ such that $\uD \psi= \frac n2 H\psi$.
\end{cor}
{\bf Proof.} To apply Theorem \ref{thm10} it suffices to prove that 
$$\mathrm{R}^\cZ + 2i\gamma(\Omega^\cZ) + 2i \gamma(N) \gamma(N\lrcorner\Omega^\cZ) $$ is non-negative. This is guarantee by the second point of Remark \ref{remarkk}. 
Proposition \ref{barconf} with the fact that $\lambda_1(\uD) = \frac n2$
imply that the immersion attains the maximum for $\lambda_1(\uD)$, since
small deformations of hypersurfaces preserve the positivity of the mean
curvature.
\begin{remark}
Note that the necessary condition $\nabla\psi =0$ obtained in
Corollary \ref{cor11} does not mean necessarily that the eigenspinor $\psi$  is the
restriction to $M$ of a parallel spinor field on $\cZ$. The reason is that  
is defined only on the hypersurface $M$. Furthermore, it is interesting
to point out that, $\nabla\psi=0$ is equivalent to 
$$\un_X\psi = -\frac 12 \ug (AX)\psi,$$
for all $X\in \Gamma(TM)$.
A spinor field $\psi$  satisfying such an overdetermined system is usually
called a generalized Killing spinor.
\end{remark}
                         
 Now, we are ready to prove Theorem \ref {thmmain2}:

{\bf Proof of Theorem \ref{thmmain2}:} We recall that the cone $\cZ$ over $P$ is the Riemannian manifold 
$\cZ= (0, \infty) \times P$  equipped with the metric
$$\<\cdot, \cdot\> = dr\otimes dr + r^2 \<\cdot, \cdot\> ,$$
for any $r \in (0, +\infty) $. In this case, we know that the vector field $x = r\frac{\partial}{\partial r} =r \partial_r$ satisfies
 $\nabla_Xx= X$ for any $X \in \Gamma(T\cZ)$. The existence of a real
 Killing spinor on $P$ implies that $\cZ$ is a $\Spinc$ manifold carrying a parallel spinor $\Psi$. Moreover, the curvature of the auxiliary line bundle associated with the 
 $\Spinc$ structure on the cone $\cZ$ satisfies \cite{IEspinc}
$$\frac {\partial}{\partial r} \lrcorner \Omega^\cZ =N \lrcorner \Omega^\cZ=0.$$

From Proposition \ref{byfirsteigenspinors}, we
know that the restriction of $\Psi$ to $M$ satifies $\uD \Psi = \frac n2 H\Psi$. Define, as in \cite{HM2012}, the sppinor $\psi$ by $\psi := \gamma(Hx+N)\Psi$. Then, for each $X\in \Gamma(TM)$, we have 
\begin{eqnarray}\label{position}
 \nabla_X\psi = \gamma(HX-AX)\Psi,
\end{eqnarray}
since H is constant. Then, if $\{e_1,\cdots, e_n\}$ is a local orthonormal basis
tangent to $M$, we have
$$\sum_{j=1}^n \gamma(e_j)\nabla_{e_j}\psi = \sum_{j=1}^n \gamma(e_j) \gamma(He_j-Ae_j)\Psi = -\tr (H\id-A)\Psi = 0.$$   
Hence, the spinor fields $\gamma(Hx+N)\Psi$, where $\Psi$ is a parallel spinor field is also an eigenspinor
corresponding to
the eigenvalue $\frac n2 H$. By Corollary \ref{cor11}, we have $\nabla (\gamma(Hx+N)\Psi)= 0$. Identity (\ref{position}) becomes  $\gamma(HX-AX)\Psi=0$ and so $AX=H X$, for all $X \in \Gamma(TM)$. This means that $M$ is a totally umbilical hypersuface of $\cZ$ and of constant mean curvature $H$.

We recall that for totally umbilical hypersurfaces of constant mean curvature of a $\Spinc$ manifold with parallel spinor, we have $\Vert \mathrm{Ric}^\cZ (N, N) \Vert = \Vert  N \lrcorner \Omega^{\cZ}\Vert$ (See Proposition \ref{hyper_totally_umb}) but here $N \lrcorner \Omega^{\cZ} = 0$  so the Ricci condition (\ref{CO}) in Theorem \ref{Montiel-classification} is satisfied ($V$ in this case is $\frac{\partial}{\partial t}$ and hence $M$ is a slice $\{ s \} \times P$.

Now, we move to prove Theorem \ref{thmmain3}:

{\bf Proof of Theorem \ref{thmmain3}:} Since $\cZ$ satisfies $\mathrm{R}^\cZ \geq c_{n+1}\vert\Omega^\cZ\vert  $
and $M$ is a compact  embedded hypersurface in $\cZ$, then the $\Spinc$ lower bound of Hijazi-Montiel-Zhang 
\cite{NR1, HMZ1} is satisfied, i.e., 
$$\lambda_1(\uD) \geq \frac n2 \inf H = \frac n2 H.$$
Since the mean curvature $H$ is constant and $\cZ$ has a parallel spinor, we get that
$$\lambda_1(\uD) \leqslant \frac n2  H.$$
Hence, $\lambda_1(\uD)  = \frac n2 H$ and M is immersed into $\cZ$ by first eigenspinors. By Theorem \ref{thmmain2}, we get the result.

\section{Non-Spin examples}\label{nonspin}
In this section, we present two low-dimensional examples of Theorem \ref{thmmain2} that are not $\Spin$.
\subsection{The $3$-dimensional homogeneous manifolds with $4$-dimensional isometry group $\Ekt$}
We denote a $3$-dimensional homogeneous manifolds with $4$-dimensional isometry group by $\Ekt$. It is a Riemannian fibration over a simply connected 2-dimensional manifold $\M^2(\kappa)$ with constant curvature $\kappa$ and such that the fibers are geodesic. We denote by $\tau$ the bundle curvature, which measures the deviation of the fibration to be a Riemannian product. Precisely, we denote by $\xi$ a unit vertical vector field, that is tangent to the fibers. The vector field $\xi$ is a Killing field and satisfies for all  vector field $X$, 
$$\nabla_X\xi=\tau X\wedge\xi,$$ 
where $\nabla$ is the Levi-Civita connection and $\wedge$ is the exterior product. When $\tau$ vanishes, we get a product manifold $\M^2(\kappa)\times\mathbb R$. If $\tau\neq0$, these manifolds are of three types: They have the isometry group of the Berger spheres if $\kappa>0$, of the Heisenberg group $\mathrm{Nil}_3$ if $\kappa=0$ or of $\widetilde{\mathrm{PSL}_2(\mathbb R)}$ if $\kappa<0$.\\\\
Note that if $\tau =0$, then $\xi =\frac{\partial}{\partial t}$ is the unit vector field giving the orientation of $\mathbb R$ in the product $\M^2 (\kappa) \times \mathbb R$. The manifold $\Ekt$, with $\tau \neq 0$,  admits a local direct orthonormal frame $\{e_1, e_2, e_3\}$ with 
$$e_3 = \xi,$$
and such that the Christoffel symbols $ \Gamma^k_{ij} = \< \nabla_{e_i}e_j, e_k\>$ are given by
\begin{eqnarray}\label{christoffel}
\left\lbrace  
\begin{array}{l}
{\Gamma}_{12}^3={\Gamma}_{23}^1=-{\Gamma}_{21}^3=-{\Gamma}_{13}^2=\tau,\\ \\
{\Gamma}_{32}^1=-{\Gamma}_{31}^2=\tau-\frac{\kappa}{2\tau}, \\ \\
{\Gamma}_{ii}^i={\Gamma}_{ij}^i={\Gamma}_{ji}^i={\Gamma}_{ii}^j=0,\quad\forall\,i,j\in\{1,2,3\},
\end{array}
\right. 
\end{eqnarray}
We call $\{e_1, e_2, e_3=\xi\}$ the canonical frame of $\Ekt$. In this canonical frame, we have
$\Ric(e_1) = (\kappa-2\tau^2)e_1, \Ric(e_2) = (\kappa-2\tau^2)e_2, \Ric(e_3) = 2\tau^2 e_3$.

We consider the $3$-dimensional Berger Sphere $\mathbb E (\kappa, \tau)= \mathbb B^3$ (because it is compact)  of scalar curvature $2\kappa-2 \tau^2$. It has a $\Spinc$ structure carrying a Kiling spinor field $\psi$ of Killing spinor $\frac {\tau}{2}$ \cite{NR, nakadthese}. Moreover, for  this $\Spinc$ structure, we have \cite{NR, nakadthese}
$$\gamma(\xi) \psi= -i  \psi,$$
$$\Omega (e_1, e_2) = -(\kappa-4 \tau^2),$$
 $$\xi \lrcorner \Omega = 0.$$
 We note that the Berger sphere does not admit a Spin structure with a real Killing spinor; otherwise, it would be an Einstein manifold. We consider now the cone over the Berger sphere given by $\cZ = (0, +\infty) \times \mathbb B^3$ with metric $\<.,.\> = dr\otimes dr + r^2 g$.  We calculate
$$\Ric^\cZ (e_1, e_1) = \Ric(e_1, e_1) -2 = \kappa-2 \tau^2 -2$$ 
$$\Ric^\cZ (e_2, e_2) = \Ric(e_2, e_2) -2 = \kappa-2 \tau^2 -2$$ 
$$\Ric^\cZ (e_3, e_3) = \Ric(e_3, e_3) -2 = 2 \tau^2 -2$$ 
$$\Ric^\cZ (\partial_r, \partial_r) = 0$$ 
In addition, 
\begin{eqnarray*}
\mathrm R^\cZ -c_4 \vert \Omega^\cZ\vert &=& \frac{1}{r^2} (\mathrm R^M -6 -2\sqrt 2 \vert \Omega^M\vert) \\ &=& \frac {1}{r^2}\big(2 \kappa-2 \tau^2 -6 -2\sqrt 2 (\kappa-4 \tau^2)\big)
\end{eqnarray*}
It is clear that Condition (\ref{co1}) and Condition  (\ref{co2}) are both satisfied if we assume that 
$$ 4\tau^{2}  < \kappa \leq (7 + 3\sqrt{2})\tau^{2} - (3 + 3\sqrt{2}),$$
for $\tau \geq 2 $.
Thus, the only compact constant mean curvature hypersurfaces of $\cZ = (0, +\infty) \times \mathbb \Ekt$ immersed by first eigenspinors,  and for $\tau \geq 2$, and 
$$4\tau^{2} < \kappa \leq (7 + 3\sqrt{2})\tau^{2} - (3 + 3\sqrt{2}),$$
are slices $\{s\} \times \mathbb B^3$.
\subsection{An explicit non-Einstein and non-Spin Sasakian $5$-dimensional manifold}
We explicitly construct a 5-dimensional \emph{non-Einstein, non-Spin Sasakian manifold} $P$ with a canonical $\Spinc$ structure carrying a Killing spinor, such that its cone $\cZ$ fulfills the conditions of Theorem \ref{thmmain2}. 

Consider the K\"ahler manifold $\mathbb{CP}^2$ with the Fubini-Study metric $\omega_{FS}$. We equip $\mathbb {CP}^2$ with the rescaled Fubini-Study metric $\tilde \omega$ given by
\[
\tilde \omega = \frac{1}{6} \, \omega_{FS},
\] 
which increases the transverse scalar curvature to $72$. Over $\mathbb {CP}^2$, we take a regular circle bundle $\pi: P \to \mathbb {CP}^2$ with very small Euler class 
\[
c_1(L) = \frac{1}{24} [\tilde \omega],
\] 
ensuring that $P$ is \emph{not Spin}. The manifold $P$ is the total space of the principal $\mathbb S^1$-bundle over $\mathbb {CP}^2$ with Euler class $c_1(L)$ . It is a $5$-dimensional closed manifold. The
 standard Sasakian metric $g$ on $P$ is non-Einstein since the Euler class does not satisfy the Sasaki-Einstein condition. It also has a canonical $\Spinc$ structure with a Killing spinor, and the curvature of the auxiliary line bundle is given by   $-i\rho^T$, where $\rho^T$ is the transverse Ricci form. Consider the cone $\cZ$ over $P$. The cone metric is given by  \[
\<.,.\> = dr\otimes dr + r^2 g.
\] 
 The Ricci tensor on $\cZ$ in an orthonormal frame $\{e_1, e_2, e_3, e_4, \xi, \partial r\}$ has components
\[
\operatorname{Ric}^{\cZ}(\partial_r,\partial_r) = 0, \quad
\operatorname{Ric}^{\cZ} (\xi,\xi) = 0, \quad
\operatorname{Ric}^{\cZ}(X,X) = 30 \quad \text{for } X=e_1, e_2, e_3 \ \text{or}\  e_4,
\] 
 The scalar curvature of the cone is then
\[
\mathrm{R}^{\cZ} =  \frac{120}{r^2}.
\]
Now we calculate
\[
\mathrm{R}^{\cZ} - 2 \sqrt{3}\, |\Omega^{\cZ}| = \frac{120 - 3\sqrt{3}}{r^2} \geq 0
\] 
It is clear that Condition (\ref{co1}) and Condition  (\ref{co2}) are both satisfied. Thus, the only compact constant mean curvature hypersurfaces of $\cZ = (0, +\infty) \times P$ immersed by first eigenspinors
are slices $\{s\} \times  P$.
                                                                                          
\end{document}